\documentclass[11pt]{amsart}
\usepackage{mathrsfs}

\usepackage{graphicx}

\usepackage{amssymb}
\usepackage{amscd}

\newtheorem{thm}{Theorem}[section]
\newtheorem{lem}[thm]{Lemma}
\newtheorem{prop}[thm]{Proposition}
\newtheorem{cor}[thm]{Corollary}

\newcommand{\bc}{\mathbb{C}}

\newcommand{\bz}{\mathbb{Z}}
\newcommand{\bp}{\mathbb{P}}

\newcommand{\co}{\mathcal{O}}

\newcommand{\cg}{\mathcal{G}}
\newcommand{\ce}{\mathcal{E}}

\newcommand{\db}{\bar\partial}

\newcommand{\cu}{L\mathfrak{U}}

\newcommand{\fh}{\mathfrak{h}}
\newcommand{\fH}{\mathfrak{H}}

\newcommand{\fu}{L\mathfrak{U}}
\newcommand{\ff}{\mathfrak{F}}

\newcommand{\pic}{\mathrm{Pic}(L\mathbb{P}_1)}

\newcommand{\picz}{\mathrm{Pic}^0(L\mathbb{P}_1)}

\newcommand{\dg}{H^{0,1}(L\bp_1)}

\newcommand{\homo}{\mathrm{Hom}(L\bc^{\ast}, \bc^{\ast})}
\newcommand{\homoc}{\mathrm{Hom}(L\bc^{\ast}, \bc)}
\newcommand{\homoz}{\mathrm{Hom}^0(L\bc^{\ast}, \bc^{\ast})}

\newcommand{\fe}{\mathfrak{E}}

\begin{document}

\title{The Picard group of the loop space of the Riemann
sphere}

\author{Ning Zhang}

\address{Department of Mathematics\\Syracuse University\\
Syracuse, NY 13244\\
U. S. A.}

\address{\emph{Present address:}
School of Mathematical Sciences\\
Peking University\\
Beijing, 100871\\ P. R. China}

\email{nzhang@math.pku.edu.cn}



\begin{abstract}
The loop space $L\bp_1$ of the Riemann sphere consisting of all
$C^k$ or Sobolev $W^{k,p}$ maps $S^1 \to \bp_1$ is an infinite
dimensional complex manifold. We compute the Picard group $\pic$ of
holomorphic line bundles on $L\bp_1$ as an infinite dimensional
complex Lie group with Lie algebra the Dolbeault group $\dg$. The
group $G$ of M\"{o}bius transformations and its loop group $LG$ act
on $L\bp_1$. We prove that an element of $\pic$ is $LG$-fixed if it
is $G$-fixed, thus completely answer the question of Millson and
Zombro about the $G$-equivariant projective embedding of $L\bp_1$.

\end{abstract}

\thanks{This work was partially supported by National
Natural Science Foundation of China grant
10871002.}

\maketitle
\section{Introduction}

Let $M$ be a finite dimensional complex manifold. We fix a
smoothness class $C^k$, $k=1,2, \cdots, \infty$, or Sobolev
$W^{k,p}$, $k=1,2, \cdots$, $1 \le p<\infty$, and consider the loop
space $LM=L_k M$, or $L_{k, p} M$, of all maps $S^1 \to M$ with the
given regularity. It is an infinite dimensional complex
Banach/Fr\'echet manifold, see \cite{l2}. The goal of this paper is
to study the Picard group $\pic$ of holomorphic line bundles on the
loop space $L\bp_1$ of the Riemann sphere.

Let $G \simeq PGL(2,\bc)$ be the group of holomorphic automorphisms
of $\bp_1$. Its loop space $LG$ with pointwise group operation is
again a complex Lie group and acts on $L\bp_1$ holomorphically, thus
also acts on $\pic$ and the Dolbeault groups of $L\bp_1$. In
\cite{mz}, Millson and Zombro raised the following question, which
is a direct motivation to study holomorphic line bundles on
$L\bp_1$: does there exist a $G$-equivariant holomorphic embedding
of $L\bp_1$ into a projectivized Banach/Fr\'echet space? In
\cite{z}, the infinite dimensional subgroup of $\pic$ of $LG$-fixed
elements was explicitly constructed; furthermore, it was proved that
the space of holomorphic sections of any line bundle in this
subgroup is finite dimensional. The Dolbeault group $\dg$ was
computed and its irreducible $G$-submodules were identified in
\cite{lz}. Based on previous works, this paper offers a complete
answer to the question of Millson and Zombro.

The following are the main results of this paper.

Let $\picz \subset \pic$ be the subgroup of topologically trivial
line bundles.

\begin{thm} \label{main0}
 There is an exact sequence of
homomorphisms
\begin{equation} \label{seq}
0 \to H^1(L\bp_1, \bz) \to \dg \stackrel{\fe}{\to} \picz \to 0.
\end{equation}
\end{thm}

The map $\fe$ in (\ref{seq}) is explicitly constructed, and it is
equivariant with respect to the group of holomorphic automorphisms
of $L\bp_1$.

Theorem \ref{main0} is closely related to the Dolbeault isomorphism.
Let $\co$ (resp. $\co^{\ast}$) be the sheaf of germs of holomorphic
(resp. non-vanishing holomorphic) functions. The short exact
sequence of sheaves
$$0 \to \bz \to \co \stackrel{\exp(2\pi i\cdot)}{\longrightarrow}\co^{\ast} \to 0$$
on a complex manifold $N$ induces a long exact sequence of sheaf
cohomology groups $$\cdots \to H^1(N, \bz) \to H^1(N,\co) \to H^1(N,
\co^{\ast}) \to H^2(N, \bz) \to \cdots,$$ where $H^1(N,\co^{\ast})
\simeq \mathrm{Pic}(N)$. Thus a map $H^{0,1}(N) \to \mathrm{Pic}(N)$
similar to $\fe$ in (\ref{seq}) would follow from the Dolbeault
isomorphism $H^{0,1}(N) \simeq H^1(N, \co)$. Such an isomorphism on
the space $L_{\infty} \bp_1$ of $C^{\infty}$ loops can be obtained
from the exactness of the sequence of sheaves $\ce_{0, 0}
\stackrel{\db}{\rightarrow} \ce_{0, 1} \stackrel{\db}{\rightarrow}
\ce_{0, 2}$, where $\ce_{p, q}$ is the sheaf of germs of
$C^{\infty}$ forms of type $(p, q)$ (see \cite[Appendix 3]{m}), and
from the existence of $C^{\infty}$ partitions of unity. However, the
above isomorphism in general fails on an infinite dimensional
complex Banach manifold or even an open subset of a complex Banach
space, see \cite{p2}, and it is not available on a general loop
space. As a consequence of Theorem \ref{main0} we obtain the
Dolbeault isomorphism $\dg \simeq H^1(L\bp_1, \co)$.

Recall that $\dg$, equipped with a natural topology, is a complex
locally convex topological space, see \cite[Section 0]{lz}.  There
is a unique complex Lie group structure on $\pic$ such that the map
$\fe$ in (\ref{seq}) is holomorphic, actually it can be considered
as $\exp(2\pi i\cdot)$, where $\exp$ is the exponential map of the
Lie group.

If a group $\cg$ acts on a set $V$, we write $V^{\cg}$ for the
$\cg$-fixed subset.

\begin{thm} \label{main}
$\pic^{G}=\pic^{LG}$; $\dg^{G}=\dg^{LG}$.
\end{thm}

By Theorem \ref{main} and \cite[Theorem 1.2]{z}, we obtain the
following

\begin{cor}
If $\Lambda \in \pic^{G}$, then $\dim H^0(L\bp_1, \Lambda) <\infty$.
\end{cor}

Therefore there does not exist a $G$-equivariant holomorphic
embedding of $L\bp_1$ into a projectivized Banach/Fr\'echet space,
otherwise the pull back of the hyperplane section bundle would be in
$\pic^{G}$, and its space of holomorphic sections would be infinite
dimensional, a contradiction.

This paper is organized as follows. In Section \ref{pre} we mainly
recall some relevant facts about loop spaces and in particular,
$L\bp_1$. It was proved in \cite{lz} that $\dg$ is isomorphic to a
space $\fH$ of normalized additive \v{C}ech 1-cocycles with respect
to a fixed open covering of $L\bp_1$. In Section \ref{map} we first
show that $\picz$ is isomorphic to a space $\fH^{\ast}$ of
normalized multiplicative \v{C}ech 1-cocycles with respect to the
same open covering, then we study the relation between $\fH$ and
$\fH^{\ast}$ and prove Theorem \ref{main0}, followed by a couple of
corollaries. At the end of this section, we briefly discuss the
naturality and some special properties of the topology of $\dg$ and
of $\pic$. Finally in Section \ref{fixed}, we study the structure of
$\pic^G$ and of $\dg^G$ and prove Theorem \ref{main}. Here a few key
objects introduced in \cite{lz} to understand $\dg$ take very
special forms and have alternative explanations when restricted to
$\dg^G$.

The author would like to thank L. Lempert and the referee for their
very helpful comments on the manuscript.

\section{Preliminaries \label{pre}}

If $\phi: M \to M'$ is a holomorphic map between finite dimensional
complex manifolds, then $L\phi: LM \ni x \mapsto \phi \circ x \in
LM'$ is holomorphic, and $L$ is functorial. Given $t \in S^1$, the
evaluation map $E_t: LM \ni x \mapsto x(t) \in M$ is holomorphic,
see \cite{l2}. The constant loops form a submanifold of $LM$, which
can be identified with $M$. If $M$ is paracompact, so is $LM$, see
\cite[Proposition 42.3]{km} for $L_{\infty}M$ and \cite[Theorem
14.17]{p} for all other cases.

Let $G \simeq PGL(2,\bc)$ be as in Section 1. If we apply the
functor $L$ to the holomorphic action $G \times \bp_1 \to \bp_1$,
then we obtain a holomorphic action $LG \times L\bp_1 \to L\bp_1$,
which induces (right) $LG$-actions on $\pic$ and the Dolbeault
groups of $L\bp_1$ by pull-backs. In particular, $G \subset LG$ acts
on $L\bp_1$.

The group $\dg$ was computed in \cite{lz} as follows. Let $\ff$ be
the space of holomorphic functions $F: \bc \times L\bc \to \bc$ with
properties
\begin{equation} \label{ff}
\left\{
\begin{array}{l}
F(\zeta/\lambda, \lambda^2y)=O(\lambda^2) \hspace{1mm} \mathrm{as}
\hspace{1mm} \bc \ni \lambda \to 0, \\
F(\zeta, x+y)=F(\zeta, x)+F(\zeta, y) \hspace{1.5mm} \mathrm{if}
\hspace{1mm} \mathrm{supp} \, x \cap \mathrm{supp} \, y =\emptyset, \\
F(\zeta, y+\mathrm{const})=F(\zeta, y).
\end{array}
\right.
\end{equation}
With the compact-open topology $\ff$ is a complex locally convex
space. There is a holomorphic (right) $G$-action on $\ff$: if $g \in
G$, write $J_g(\zeta)=d(g\zeta)/d\zeta$, and set
\begin{equation} \label{gaction}
(Fg)(\zeta, y)=F(g\zeta, y/J_g(\zeta))J_g(\zeta)
\end{equation}
(which is defined for $g\zeta\not=\infty$ and can be extended to
$\bc \times L\bc$ by the first property of $F$). Then there is a
$G$-isomorphism $\dg \simeq \bc \oplus \ff$ of locally convex
spaces, where $G$ acts on the right by the second variable.

If $a, b, \cdots \in \bp_1$, set $U_{ab\cdots}=\bp_1 \setminus \{a,
b, \cdots\}$. Thus
\begin{equation} \label{cover}
\fu=\{LU_a: a \in \bp_1\}
\end{equation}
is an open covering of $L\bp_1$. We consider $\bp_1$ as $\bc \cup
\{\infty\}$ and identify $LU_{\infty}$ with $L\bc$, a Fr\'echet
algebra. If $g \in G$, then $g(LU_{ab\cdots})=LU_{g(a)g(b)\cdots}$.

Given $v: \bp_1 \to \bc$, finitely many $a, b, \cdots \in \bp_1$ and
a function $u: LU_{ab\cdots} \to \bc$, we say that $u$ is cuspidal
with respect to $v$ (as formulated in \cite[Section 3]{lz}) if
\begin{equation} \label{cuspidal}
\lim_{\bc \ni \lambda \to \infty}
u\big(g(x+\lambda)\big)=v\big(g(\infty)\big)
\end{equation}
for any $x \in L\bc$ and $g \in G$ that maps $\infty$ to one of $a,
b, \cdots$. If $u$ is holomorphic, then it follows from
(\ref{cuspidal}), where we set $x=0$, that $u|_{U_{ab...}}$ extends
to be a holomorphic function on $\bp_1$ which agrees with $v$ at $a,
b, \cdots$, thus $u|_{U_{ab...}} \equiv c$ for some constant $c$ and
$u$ is cuspidal with respect to the constant $c$. Let
\begin{equation} \label{ind}
\mathrm{ind}_{ab}: LU_{ab} \ni x \mapsto \mathrm{ind}_{ab}(x) \in
\bz
\end{equation} (the winding number of $x: S^1 \to U_{ab}$). Then
$\mathrm{ind}_{ab}$ is cuspidal with respect to the constant $0$,
and $(\mathrm{ind}_{ab})_{a, \,b \in \bp_1}$ is a $G$-fixed \v{C}ech
1-cocycle with respect to the covering $\fu$ as in (\ref{cover}).

\begin{prop} \label{topo}
The group $H^1(L\bp_1, \bz) \simeq \bz$ is generated by the \v{C}ech
cohomology class of $(\mathrm{ind}_{ab})$. Suppose $e: \bp_1 \to
L\bp_1$ is a holomorphic embedding and there exists $t \in S^1$ such
that $E_t \circ e$, where $E_t$ is the evaluation, is an isomorphism
of $\bp_1$. Then $e^{\ast}: H^2(L\bp_1, \bz) \to H^2(\bp_1, \bz)$ is
an isomorphism.
\end{prop}
\begin{proof}
Recall that $H^1(L\bp_1, \bz)\simeq \bz$ and $H^2(L\bp_1,
\bz)\simeq\bz$ (e.g. see \cite[Part II, Proposition 15.33]{cj} and
\cite[Theorem 13.14]{p}). As every element of the open covering
$\cu$ is contractible, we have $H^1(\cu, \bz)$ $\simeq$ $H^1(L\bp_1,
\bz)$. Let $[\rho]$ be a generator of $H^1(\cu, \bz)$, where $\rho$
is a 1-cocycle. We can write $(\mathrm{ind}_{ab})=m\rho+\rho'$,
where $m \in \bz$ and $\rho'$ is exact. Since $\mathrm{ind}_{\infty
0}$ takes value of any integer and $\rho'_{\infty 0}$ is constant,
we must have $|m|=1$. It is clear that the map $e^{\ast}$ is
surjective, and any surjective homomorphism $\bz \to \bz$ is an
isomorphism.
\end{proof}

\begin{cor} \label{trivial bundle}
Let $e$ be as in Proposition \ref{topo}. Then a bundle $\Lambda \in
\pic$ is in $\picz$ if and only if $\Lambda|_{e(\bp_1)}$ (or in
particular $\Lambda|_{\bp_1}$) is trivial.
\end{cor}
\begin{proof}
The conclusion immediately follows from Proposition \ref{topo} and
the relation of the first Chern classes $c_1(e^{\ast}
\Lambda)=e^{\ast} c_1(\Lambda)$.
\end{proof}

The proof of the following proposition is similar to that of
\cite[Proposition 2.3]{z}.

\begin{prop} \label{patch}
Let $M \stackrel{\pi}{\to} B$ be a holomorphic fiber bundle with
fibers $F_b$ finite dimensional connected compact complex manifolds,
and $E \to M$ a holomorphic line bundle. If $\pi$ has a holomorphic
section $\sigma: B \to M$ such that $E|_{\sigma(B)}$ and $E|_{F_b}$,
for all $b \in B$, are holomorphically trivial, then $E$ is
holomorphically trivial.
\end{prop}
\begin{proof}
Let $u$ be a non-vanishing holomorphic section of $E|_{\sigma(B)}$,
and let $s$ be the section of $E \to M$ such that $s|_{F_b}$ is the
holomorphic section of $E|_{F_b}$ satisfying
$s(\sigma(b))=u(\sigma(b))$. It is clear that $s$ is non-vanishing.
Next we show that $s$ is $C^{\infty}$ and in turn holomorphic.

Let $b_0 \in B$. By \cite[Proposition 5.1]{l1} (where we choose
$p=0$, $q=1$ and $f=0$), there exist a neighborhood $b_0 \in B_0
\subset B$ and a section $v \in C^{\infty}(\pi^{-1}(B_0), E)$ such
that $v|_{F_b}$ is holomorphic for any $b \in B_0$, and
$v|_{F_{b_0}} \not= 0$. By choosing a sufficiently small $B_0$ we
can assume that $\pi^{-1}(B_0) \stackrel{\pi}{\to} B_0$ is trivial
and $v$ is non-vanishing. The function $s/v$ on $\pi^{-1}(B_0)$ is
constant on each fiber $F_b$, and can also be considered as a
function on $B_0$. Since $s/v(b)$ $=$ $u(\sigma(b))/v(\sigma(b))$,
it follows that $s/v$ is $C^{\infty}$. Therefore $s$ is $C^{\infty}$
and holomorphic on each fiber $F_b$ as well as the cross section
$\sigma(B)$. Apply \cite[Proposition 5.2 (ii)]{l1} to the (0,1) form
$\db s$, we obtain that $\db s=0$.
\end{proof}

\section{The map $\fe: H^{0,1}(L\bp_1) \to \mathrm{Pic^0}(L\bp_1)$
\label{map}}

In this section we construct the map $\fe$ in (\ref{seq}) and prove
Theorem \ref{main0}, then we compute $H^1(L\bp_1, \co)$ and $\pic$
as corollaries. Finally we study certain special features of the
topology of $\dg$ and of $\pic$.

Let $\fH$ (resp. $\fH^{\ast}$) be the linear space (resp. group) of
holomorphic additive (resp. multiplicative) \v{C}ech 1-cocycles
$(\fh_{ab} \in \co(LU_{ab}))_{ a, \,b \in \bp_1} $ (resp.
$(\fh^{\ast}_{ab} \in \co^{\ast}(LU_{ab}))_{ a, \,b \in \bp_1}$)
 with respect to the covering $\fu$ of $L\bp_1$
 such that every component $\fh_{ab}$ (resp.
 $\fh^{\ast}_{ab}$) is
 cuspidal with respect to $v\equiv 0$ (resp. $v\equiv 1$), see
  (\ref{cuspidal}). The cocycle
 $(\mathrm{ind}_{ab}) \in \fH$. The cuspidal
 property implies that $\fh_{ab}|_{U_{ab}} \equiv 0$ (resp.
 $\fh^{\ast}_{ab}|_{U_{ab}} \equiv 1$). The group $G$ acts on
 $\fH$ (resp. $\fH^{\ast}$) by pull-backs.

 The space $\fH$ was defined in \cite[Section 3]{lz} with the
 additional requirement that $\fh_{ab}(x)$ is holomorphic in $a$ and
 $b$,  but we shall show in Corollary \ref{holo} that it is redundant.

It was proved in \cite{lz} that $\dg \simeq \fH$ as $G$-modules. We
shall construct a $G$-isomorphism $\picz$ $\simeq$ $\fH^{\ast}$, and
relate $\dg$ and $\picz$ by the exponential map $\fH \to
\fH^{\ast}$.

We begin with a family of normalized non-vanishing holomorphic
sections of a given bundle $\Lambda \in \picz$ over $LU_a$, $a \in
\bp_1$, which depend on $a$ holomorphically.

\begin{prop} \label{section}
Let $\Lambda \in \picz$, $v$ a holomorphic section of
$\Lambda|_{\bp_1}$ and $x_0 \in L\bc^{\ast}$ (where $\bc^{\ast}=\bc
\setminus \{0\}$). For any $a \in \bp_1$ there is a unique
$\sigma_a$ $=$ $\sigma_{a,\,v,\,x_0} \in H^0(LU_a, \Lambda)$ such
that
\begin{equation} \label{lim}
\lim_{\bc \ni \lambda \to \infty} \sigma_a \left(g(x+\lambda
x_0)\right)=v(a)
\end{equation}
for all $x \in LU_{\infty}$ and $g \in G$ with $g(\infty)=a$. If $v$
is non-vanishing, so is $\sigma_a$; and $\sigma_{a}$ depends
linearly on $v$. Furthermore, $\sigma(a,y)=\sigma_a(y)$ is
holomorphic in $(a,y)$.
\end{prop}
\begin{proof}
Let $\phi$ be a continuous complex linear functional on
$LU_{\infty}$ such that $\phi(x_0)=1$, let $Z=\ker(\phi)$ and
consider the holomorphic map
\begin{equation} \label{p}
P: \bp_1 \times G \times Z \ni (\lambda, g, z) \mapsto g(z+\lambda
x_0) \in L\bp_1
\end{equation}
and the pull-back bundle $P^{\ast}\Lambda$. Note that for any $g \in
G$ and $z \in Z$, the map $P|_{\bp_1 \times \{g\} \times \{z\}}$ is
an embedding with image curve $C(g, z)$ passing through $g(\infty)
\in \bp_1 \subset L\bp_1$. The topologically trivial line bundle
$\Lambda|_{C(g, z)}$ is also holomorphically trivial, so is
$P^{\ast}\Lambda|_{\bp_1 \times \{g\} \times \{z\}}$. Thus an
equivalent statement to (\ref{lim}) is that along each curve $C(g,
z)$, where $g(\infty)=a$, $\sigma_a$ extends to be the unique
holomorphic section of $\Lambda|_{C(g, z)}$ with value $v(a)$ at
$a$.

The uniqueness of $\sigma_a$ follows from the fact that, for any
fixed $g \in G$ with $g(\infty)=a$, the set $LU_a\cup \{a\}$ is
covered by curves of the form $C(g, z)$, and $\sigma_a$ is the
unique holomorphic section of $\Lambda$ on each curve satisfying
(\ref{lim}). If $v(a) \not =0$, then $\sigma_a$ is non-vanishing on
each curve, and $\sigma_a$ depends on $v(a)$ linearly on each curve.

As to existence, since $P$ maps $\{\infty\} \times G \times Z$ to
$\bp_1 \subset L\bp_1$, the bundle $P^{\ast}\Lambda$ is trivial on
$\{\infty\} \times G \times Z$, which can be considered as a cross
section of the trivial $\bp_1$-bundle $\bp_1 \times G \times Z$. It
follows from Proposition \ref{patch} that $P^{\ast} \Lambda$ is
trivial and we can uniquely extend $P^{\ast} v$ to a section
$\tilde{\sigma} \in H^0(\bp_1 \times G \times Z, P^{\ast}\Lambda)$.
Note that, for any fixed $g \in G$ with $g(\infty)=a$, $P$ maps $\bc
\times \{g\} \times Z$
 biholomorphically onto $LU_a$. Let
$\tilde{P}: P^{\ast}\Lambda \to \Lambda$ be the bundle map
associated with $P$. Then there exists $\sigma_a \in H^0(LU_a,
\Lambda)$ such that
 \begin{equation} \label{sigma} \sigma_a \circ P=
 \tilde{P} \circ \tilde{\sigma} \end{equation}
on $\bc \times \{g\} \times
 Z$. It follows from the property of $\tilde{\sigma}$ that
 $\sigma_a$ satisfies (\ref{lim}).

Choose $g \in G$ with $g(\infty)=a$
 such that $g$ depends holomorphically on $a$ (which can be done
 locally). If $y=P(\lambda, g, z) \in LU_a$, then we can write
 $\lambda$ and $z$ in terms of $y$ and $g$, and equation (\ref{sigma})
 becomes
$$\sigma_a(y)=\tilde{P}\left(\tilde{\sigma}
 \left(\phi\left(g^{-1}(y)\right), g,
 g^{-1}(y)-\phi\left(g^{-1}(y)\right)x_0\right)\right),$$
 i.e. $\sigma_a(y)$ is
 holomorphic in $(a,y)$.
\end{proof}

When $\Lambda$ is $LG$-fixed, two sections $\sigma_{a, \,v, \,x_0}$
and $\sigma_{a, \,v, \, x_1}$, $x_0, x_1 \in L\bc^{\ast}$, as in
Proposition \ref{section} only differ by a multiplicative constant,
see \cite[Proposition 2.7]{z}. However, this is not true for a
general $\Lambda \in \picz$.

Choose $v$ to be non-vanishing and set $x_0=1$ in Proposition
\ref{section}. Since $\Lambda$ determines $v$ up to a multiplicative
constant, we can uniquely associate with $\Lambda$ the \v{C}ech
1-cocycle $(\fh^{\ast}_{ab}=\sigma_b/\sigma_a)$. The property
(\ref{lim}) of $\sigma_a$ and $\sigma_b$ implies that
$(\fh^{\ast}_{ab}) \in \fH^{\ast}$.

\begin{prop} \label{iso}
The $G$-morphism $\picz \ni \Lambda \mapsto (\fh^{\ast}_{ab}) \in
\fH^{\ast}$ is an isomorphism of groups.
\end{prop}
\begin{proof}
It is clear that the kernel only contains the trivial bundle. Given
$(\fh^{\ast}_{ab}) \in \fH^{\ast}$, we can construct a line bundle
$\Lambda$ by taking the union of $LU_a \times \bc$ over all $a \in
\bp_1$ and identifying $\{x\} \times \bc$ in $LU_a \times \bc$ and
$LU_b \times \bc$ via multiplication by $\fh^{\ast}_{ab}(x)$. Let
$\sigma_a$ be the section of $\Lambda$ on $LU_a$ corresponding to
$LU_a \times \{1\}$. Then $\fh^{\ast}_{ab}=\sigma_b/\sigma_a$. Since
$\fh^{\ast}_{ab}|_{U_{ab}} \equiv 1$, there is a unique
non-vanishing holomorphic section $v$ of $\Lambda|_{\bp_1}$ such
that $v|_{U_a}=\sigma_a|_{U_a}$. It follows from Corollary
\ref{trivial bundle} that $\Lambda \in \picz$. From the relation
$\sigma_a =\sigma_b/\fh^{\ast}_{ab}$ and the cuspidal property of
$\fh^{\ast}_{ab}$ we can obtain (\ref{lim}).
\end{proof}

Any $(\fh_{ab}) \in \fH$ (resp. $(\fh^{\ast}_{ab}) \in \fH^{\ast}$)
can be considered as a function on
$$\Omega=\{(a,b,x) \in \bp_1 \times \bp_1 \times L\bp_1: a, b
\not\in x(S^1)\}.$$

\begin{cor} \label{holo}
If $(\fh_{ab}) \in \fH$ (resp. $(\fh^{\ast}_{ab}) \in \fH^{\ast}$),
then $\fh_{ab}(x)$ (resp. $\fh^{\ast}_{ab}(x)$) is holomorphic in
$(a, b, x) \in \Omega$.
\end{cor}
\begin{proof}
By Propositions \ref{iso} and \ref{section} there exist $\Lambda \in
\picz$ and sections $\sigma_a \in H^0(LU_a, \Lambda)$, $a \in
\bp_1$, such that $\fh^{\ast}_{ab}=\sigma_b/\sigma_a$, and
$\sigma_a(x)$ is holomorphic in $(a, x)$. The conclusion immediately
follows for $(\fh^{\ast}_{ab})$. It is also true for $(\fh_{ab})$
because $(\fh_{ab})$ is a logarithm of $(\exp(\fh_{ab})) \in
\fH^{\ast}$.
\end{proof}

As a subspace of $\co(\Omega)$ with the compact-open topology, $\fH$
is a complete locally convex space. The isomorphism $\dg \simeq \fH$
is topological, see \cite[Section 3]{lz}.

Next we study the exponential map $\fH \to \fH^{\ast}$.

\begin{prop} \label{logarithm}
Suppose that we are given finitely many $a, b, \cdots \in \bp_1$,
and a function $h^{\ast}$ $\in$ $\co^{\ast}(LU_{ab\cdots})$ cuspidal
with respect to a constant $c^{\ast}$. If there exists $h \in
\co(LU_{ab\cdots})$ such that $\exp(h)=h^{\ast}$, then $h$ is
cuspidal with respect to a constant $c$, where $\exp(c)=c^{\ast}$.
\end{prop}
\begin{proof}
The cuspidal property of $h^{\ast}$ implies that
$h^{\ast}|_{U_{ab\cdots}} \equiv c^{\ast}$. So $h|_{U_{ab\cdots}}
\equiv c$ for some constant $c$, where $\exp(c)=c^{\ast}$. Take any
point from $\{a, b, \cdots\}$, say $a$, and fix $g \in G$ with
$g(\infty)=a$. Since $h$ is a logarithm of $h^{\ast}$, the limit
$$\omega(x) \triangleq \lim_{\bc \ni \lambda \to \infty}
h(g(x+\lambda))$$ exists for any $x \in L\bc$, $\exp(\omega(x))
\equiv c^{\ast}$ and $\omega(0)=c$. Next we show that $\omega(x)
\equiv c$.

Recall the map $P$ and the hyperplane $Z \subset L\bc$ as in
(\ref{p}), where we set $x_0 =1$.  It follows from the definition of
$\omega$ that, for any $x \in L\bc$, there exists $z\in Z$ such that
$\omega(x)=\omega(z)$. We only need to show that $\omega$ is
continuous, thus constant on any line $\bc z$ spanned by $z \in Z
\setminus \{0\}$, so $\omega(z)=\omega(0)=c$.

Let $P_{g,\,z}: \bp_1 \times \bc z \to L\bp_1$ be the map obtained
from $P$ with the fixed $g$ and $z$ above. Note that $P_{g, \,z}$
maps a neighborhood $U$ of the set $\{\infty\} \times \bc z$ to
$LU_{ab\cdots} \cup \{a\}$, and a point of $U$ is mapped to $a$ if
and only if it is in $\{\infty\} \times \bc z$. Define the function
$\tau^{\ast}$ on $U$ as in the following: $\tau^{\ast}$ takes value
$c^{\ast}$ on $\{\infty\} \times \bc z$, and $\tau^{\ast}=P_{g, \,
z}^{\ast}h^{\ast}$ otherwise. It follows from the cuspidal property
of $h^{\ast}$ that $\tau^{\ast}$ is holomorphic with respect to each
variable, thus holomorphic on $U$. Given $\zeta_0 \in \bc$, we
choose a small neighborhood $V \subset U$ of $(\infty, \zeta_0 z)$
such that $V$ is simply connected, and $V \setminus \{\infty\}
\times \bc z$ is connected. Let $\tau$ be the logarithm of
$\tau^{\ast}$ on $V$ such that $\tau$ agrees with $P_{g, \,z}^{\ast}
h$ on $V \setminus \{\infty\} \times \bc z$. Then $\omega(\zeta
z)=\tau(\infty, \zeta z)$ for $\zeta$ near $\zeta_0$, and $\omega$
is continuous on $\bc z$.
\end{proof}

\begin{lem} \label{key}
The $G$-morphism $\fe: \fH \ni (\fh_{ab}) \mapsto \left(\exp(2\pi
i\fh_{ab})\right) \in \fH^{\ast}$ is surjective, and its kernel
consists of integer multiples of the cocycle $(\mathrm{ind}_{ab})$
as in (\ref{ind}).
\end{lem}
\begin{proof} Let $\fh^{\ast} \in \fH^{\ast}$. Since the line bundle
associated with $\fh^{\ast}$ is topologically trivial, there exist
non-vanishing functions $\psi^{\ast}_a \in C(LU_a)$, $a \in \bp_1$,
such that $\psi^{\ast}_b/\psi^{\ast}_a=\fh^{\ast}_{ab}$. As
$\fh^{\ast}_{ab}|_{U_{ab}} \equiv 1$, the functions
$\psi^{\ast}_a|_{U_a}$, $a \in \bp_1$, form a non-vanishing
continuous function $v$ on $\bp_1 \subset L\bp_1$. Fix a logarithm
$u$ of $v$. Since $LU_a$, $a \in \bp_1$, are simply connected, there
exist $\psi_a \in C(LU_a)$ such that $\exp(2\pi i
\psi_a)=\psi^{\ast}_a$ and $2\pi i\psi_a|_{U_a}=u|_{U_a}$. Consider
the \v{C}ech 1-cocycle $(\fh_{ab}=\psi_b-\psi_a)$. It is clear that
$\fh_{ab}|_{U_{ab}} \equiv 0$ and $\exp(2\pi
i\fh_{ab})=\fh^{\ast}_{ab}$, so $\fh_{ab}$ $\in$ $\co(LU_{ab})$. It
follows from Proposition \ref{logarithm} that $\fh_{ab}$ is cuspidal
with respect to the constant $0$, therefore $(\fh_{ab}) \in \fH$.

The cocycle $(\mathrm{ind}_{ab})$ is clearly in the kernel of $\fe$.
Suppose $\fh \in \fH$ is in the kernel of $\fe$, then it must be
$\bz$-valued, thus generates a cohomology class in $H^1(\fu, \bz)$
$\simeq$ $H^1(L\bp_1, \bz)$. It follows from Proposition \ref{topo}
that $\fh=m (\mathrm{ind}_{ab})+\fh'$, where $m \in \bz$ and $\fh'$
is an exact $\bz$-valued cochain of $\fu$, in particular, any
component $\fh'_{ab}$ is constant. The cuspidal property of
$\fh'_{ab}$ implies that $\fh'_{ab} \equiv 0$.
\end{proof}

Now $\fH^{\ast}$, the quotient space of the complex locally convex
space $\fH$ modulo a discrete subgroup, is a complex Lie group with
the quotient topology. We shall show in Proposition \ref{2t} that
this topology of $\fH^{\ast}$ is actually the compact-open topology.

Theorem \ref{main0} immediately follows from \cite[Theorem 3.3]{lz}
(that $\dg \simeq \fH$), Proposition \ref{iso}, Lemma \ref{key} and
Proposition \ref{topo}. The map $\fe$ in Theorem \ref{main0} is the
one in Lemma \ref{key}. If we endow $\picz$ with the complex Lie
group structure from $\fH^{\ast}$, then the map $\fe$ in Theorem
\ref{main0} is just $\exp(2\pi i \cdot)$, where $\exp$ is the
exponential map of the Lie group. There is a unique complex Lie
group structure on $\pic$ such that $\picz$ is the component
containing $ 1 \in \pic$.

Next we draw a couple of corollaries of Lemma \ref{key}/Theorem
\ref{main0}.

Note that any 1-cocycle in $\fH$ (resp. $\fH^{\ast}$) generates a
\v{C}ech cohomology class in $H^1(\fu, \co)$ (resp.
$H^1(\fu,\co^{\ast})$). From the definition of $\fH$ (resp.
$\fH^{\ast}$) we obtain the following

\begin{prop}
The map $\fH \to H^1(\fu, \co)$ (resp. $\fH^{\ast}$ $\to$
$H^1(\fu,\co^{\ast}))$ is injective.
\end{prop}
\begin{proof}
If $\{\fh_{ab}\} \in \fH$ is exact, i.e. there exist $h_a \in
\co(LU_a)$, $a \in \bp_1$, such that $\fh_{ab}=h_a-h_b$, it follows
from the cuspidal property of $\fh_{ab}$ that
$$\lim_{\bc \ni \lambda \to \infty} h_a\big(g(x+\lambda)\big)
=\lim_{\lambda \to
\infty}(\fh_{ab}+h_b)\big(g(x+\lambda)\big)=h_b(a),$$ for any $x \in
L\bc$ and $g \in G$ with $g(\infty)=a$. So $h_a$ is constant along
any curve of the form $\{g(x+\lambda): \lambda \in \bp_1\} \subset
L\bp_1$, and $h_a\equiv h_b(a)$. Also $h_b$ is constant, thus
$\fh_{ab}\equiv 0$. The conclusion for $\fH^{\ast}$ follows from
similar arguments.
\end{proof}

Therefore we can consider $\fH$ (resp. $\fH^{\ast}$) as a subspace
(resp. subgroup) of $H^1(\fu, \co)$ (resp. $H^1(\fu,\co^{\ast})$),
which is, in turn, a subspace (resp. subgroup) of $H^1(L\bp_1, \co)$
(resp. $H^1(L\bp_1, \co^{\ast}$)).

\begin{cor} \label{dol}
The maps $\fH \to H^1(\fu, \co) \to H^1(L\bp_1, \co)$ are
isomorphisms of groups. In particular, $\dg$ $\simeq$ $H^1(L\bp_1,
\co)$.
\end{cor}
\begin{proof}
Since $L\bp_1$ is rationally connected, see the proof of Proposition
\ref{section}, any holomorphic function on $L\bp_1$ is constant.
Therefore $H^0(L\bp_1, \co)$ $\simeq$ $\bc$ and $H^0(L\bp_1,
\co^{\ast})$ $\simeq$ $\bc^{\ast}$. The short exact sequence of
sheaves $ 0\to \bz \to \co \stackrel{\exp(2\pi
i\cdot)}{\longrightarrow} \co^{\ast} \to 0$ induces a long exact
sequence of cohomology groups
\begin{eqnarray}
& 0\to \bz \to \bc \to \bc^{\ast} \to H^1(L\bp_1, \bz) \to
H^1(L\bp_1, \co) \nonumber
\\& \stackrel{\exp(2\pi i\cdot)^{\ast}}{\longrightarrow}
H^1(L\bp_1, \co^{\ast}) \to H^2(L\bp_1, \bz) \to \cdots,
\label{longexact}
\end{eqnarray}
where the kernel of $\exp(2\pi i\cdot)^{\ast}$ is
$H^1(L\bp_1,\bz)$$\subset \fH \subset$ $H^1(L\bp_1, \co)$, and the
range of $\exp(2\pi i\cdot)^{\ast}$ is $\picz \subset H^1(L\bp_1,
\co^{\ast})$. It follows from Lemma \ref{key} and Propositions
\ref{topo} and \ref{iso} that the map $\exp(2\pi i\cdot)^{\ast}$ and
its restriction to $\fH$ have the same kernel
 and the same range, thus the map
$\fH \to H^1(L\bp_1, \co)$ is an isomorphism.
\end{proof}

With the Dolbeault isomorphism in Corollary \ref{dol}, the map $\fe$
in Theorem \ref{main0} can be considered as $\exp(2\pi
i\cdot)^{\ast}$ in (\ref{longexact}), thus it is equivariant with
respect to the group of holomorphic automorphisms of $L\bp_1$.

Now $\picz$ can be computed  as a quotient group by Theorem
\ref{main0}. To identify $\picz$ more precisely, we recall a few
concepts which played key roles in \cite{lz}, and which are also
useful in the next section. Let $\ff$ be the function space as in
(\ref{ff}). In \cite[Section 4]{lz}, a map $\label{a} \alpha: \fH
\to \ff$ was introduced as follows. Given $\fh=(\fh_{ab}) \in \fH$,
the cocycle relation implies that $d_{\zeta}\fh_{a\zeta}(x)$ is
independent of $a$. For $\zeta \in \bc$ we can write
\begin{equation} \label{f}
d_{\zeta}\fh_{a\zeta}(x)=F\left(\zeta,
\frac{1}{\zeta-x}\right)d\zeta, \hspace{2mm} x \in LU_{\zeta},
\end{equation}
where $F \in \co(\bc \times L\bc)$. It turns out that $F \in \ff$,
and we define $\alpha(\fh)=F$. The map $\alpha$ is a continuous
$G$-morphism. It was proved in \cite[Section 5]{lz} that the kernel
of $\alpha$ is one-dimensional, spanned by the 1-cocycle
$(\mathrm{ind}_{ab})$ as in (\ref{ind}), and $\alpha$ has a right
inverse $\beta$, which is continuous and $G$-equivariant.

Since the kernel of $\exp=\fe(\cdot/2\pi i):$ $\fH \to \fH^{\ast}$
is contained in the kernel of $\alpha: \fH \to \ff$, we have a
well-defined continuous $G$-morphism $\alpha^{\ast}: \fH^{\ast} \to
\ff$ such that
\begin{equation} \label{diag}
\alpha^{\ast} \circ \exp =\alpha.
\end{equation}

\begin{prop} \label{a2}
The map $\alpha^{\ast}$ has a continuous $G$-equivariant right
inverse $\beta^{\ast}$, and its kernel consists of cocycles of the
form
\begin{equation} \label{k2} \left(\zeta^{\mathrm{ind}_{ab}}\right),
 \hspace{2mm}
\zeta \in \bc^{\ast}.
\end{equation}
\end{prop}
\begin{proof}
The right inverse $\beta^{\ast}$ is the composition of the right
inverse of $\alpha$ and the map $\exp: \fH \to \fH^{\ast}$. It is
straightforward to verify that elements of
$\ker(\alpha^{\ast})=\exp(\ker(\alpha))$ take the form as in
(\ref{k2}).
\end{proof}

Since $H^2(L\bp_1, \bz) \simeq \bz$, we obtain the following

\begin{cor} \label{pic}
$\pic \simeq \bc^{\ast} \times \ff \times \bz$.
\end{cor}

Recall that $\fH$ (resp. $\fH^{\ast}$) is a subset of $\co(\Omega)$
(resp. $\co^{\ast}(\Omega)$). In the remainder of this section we
are going to study the compact-open topology on $\fH$ and on
$\fH^{\ast}$. If we consider the 1-cocycle $(\mathrm{ind}_{ab})$ as
a function $\mathrm{ind}:\Omega \to \bz$, then for any $m \in \bz$,
the subset $\mathrm{ind}^{-1}(m) \not= \emptyset$ is a path
connected component of $\Omega$.

Let $X \subset \bc \times L\bc$ be a compact subset and let $X'$ be
the inverse image of $X$ under the homeomorphism
$$\{(\zeta, x) \in \bc \times L\bp_1: x \in LU_{\zeta} \}
 \ni (\zeta, x) \mapsto \big(\zeta, 1/(\zeta-x) \big) \in \bc \times L\bc$$
(see (\ref{f})).  If $0 < r <\inf_{(\zeta, \,x) \in \,X', \,t \in
S^1} |\zeta-x(t)|$ (where $|\zeta-x(t)|=+\infty$ when $x(t)=\infty
\in \bp_1$), then
\begin{equation}
K=K(X', r)=\{(\zeta, \zeta+z, x): (\zeta, x) \in X', z \in \bc, |z|
\le r\}
\end{equation}
is a compact subset of $\{(a, b, x) \in \Omega: a, b \not=\infty,
\mathrm{ind}_{ab}(x)=0\}$.

\begin{prop} \label{estimate}
With notation above, we have
\begin{equation} \label{ah}
\sup_{X}|\alpha(\fh)| \le r^{-1} \sup_K|\fh|
\end{equation}
for any $\fh \in \fH$; and if $0< \delta <1$, then there exists a
constant $C_2=C_2(\delta)>0$ such that
\begin{equation} \label{ash}
\sup_K \left| \fh \right| \le C_2 \sup_K \left|\exp \fh -1 \right|
\end{equation}
for any $\fh \in \fH$ with $\sup_K|\exp \fh -1| \le \delta$.
\end{prop}
\begin{proof}
Taking the derivative of $\fh(\zeta, \zeta+z, x)$ with respect to
$z$ at $z=0$ for any $(\zeta, x) \in X'$, the inequality (\ref{ah})
follows from (\ref{f}) and Cauchy's estimate. Note that $\fh (\zeta,
\zeta, x)=0$ for all $(\zeta, x) \in X'$, i.e. on every path
connected component of $K$ we can find a point where $\fh$ vanishes.
Thus (\ref{ash}) follows from basic properties of the branch of $\ln
w$ on $\{w \in \bc: |w-1|< \delta\}$ sending $w=1$ to $0$.
\end{proof}

\begin{prop} \label{2t}
Let $\omega_0 \in \Omega \setminus \mathrm{ind}^{-1}(0)$ be fixed.
Then the following topologies of $\fH$ (resp. $\fH^{\ast}$) are the
same.

\begin{enumerate}
  \item[(a)] The compact-open topology of $\fH$ (resp. the
  quotient topology of $\fH^{\ast}$ induced by the homomorphism
$\fe: \fH \to \fH^{\ast}$).
  \item[(b)] The topology that has as a basis the sets $\{\fh \in \fH:
\sup_{K \cup \{\omega_0\}} |\fh-\fh_0| < \varepsilon\}$ (resp.
$\{\fh^{\ast} \in \fH^{\ast}: \sup_{K \cup \{\omega_0\}}
|\fh^{\ast}-\fh_0^{\ast}|< \varepsilon\}$), where $K \subset \{(a,
b, x) \in \Omega: a, b \not=\infty, \mathrm{ind}_{ab}(x)=0\}$ is
compact, $\fh_0 \in \fH$ (resp. $\fh_0^{\ast} \in \fH^{\ast}$), and
$\varepsilon>0$.
\end{enumerate}
\end{prop}
\begin{proof}
Let $\mathscr{T}_a$ and $\mathscr{T}_b$ denote the topologies in (a)
and (b) respectively. It follows from Proposition \ref{estimate}
that the homomorphism $\alpha: (\fH, \mathscr{T}_b) \to \ff$ (resp.
$\alpha^{\ast}: (\fH^{\ast}, \mathscr{T}_b) \to \ff$) between
topological vector spaces (resp. topological groups) is continuous.
With the relative topology its kernel is isomorphic to $\bc$ (resp.
$\bc^{\ast}$), and its right inverse $\beta: \ff \to (\fH,
\mathscr{T}_b)$ (resp. $\beta^{\ast}: \ff \to (\fH^{\ast},
\mathscr{T}_b)$) is also continuous (for $\mathscr{T}_a$ is finer
than $\mathscr{T}_b$). So these homomorphisms induce isomorphisms
$(\fH, \mathscr{T}_b) \simeq \bc \times \ff \simeq (\fH,
\mathscr{T}_a)$ (resp. $(\fH^{\ast}, \mathscr{T}_b)$
$\simeq$$\bc^{\ast} \times \ff \simeq (\fH^{\ast}, \mathscr{T}_a)$).
\end{proof}

Since the compact-open topology of $\fH^{\ast}$ is finer than the
one in Proposition \ref{2t}(b) and coarser than the other in (a),
they are all the same.

\section{The structure of $G$-fixed subgroups \label{fixed}}

In this section we study the structure of $\dg^G$ and of $\pic^G$,
and in particular, prove Theorem \ref{main}. We begin with explicit
constructions of all elements of $\fH^G$ and $(\fH^{\ast})^G$.

Let $\homoc$ (resp. $\homo$) be the linear space (resp. group) of
all holomorphic homomorphisms from the loop group $L\bc^{\ast}$ to
$\bc$ (resp. $\bc^{\ast}$). With the compact-open topology $\homoc$
is a complex locally convex space. If $\psi \in \homoc$ and $\varphi
\in \homo$, then $\psi(z) = 0$ and $\varphi(z)=z^n$ for $z \in
\bc^{\ast} \subset L\bc^{\ast}$, where $n$ is a fixed integer, and
we call this $n$ the order of $\varphi$. We write $\homoz$ for the
subgroup of zero order elements of $\homo$.

Let $(L\bc)'$ be the space of continuous linear functionals on
$L\bc$ and let $$\mathcal{A}=\{\phi \in (L\bc)': \phi|_{\bc}\equiv
0\}.$$ Associated with any $\psi \in \homoc$ there is a commutative
diagram

\begin{equation} \label{psi}
\setlength{\unitlength}{1.1cm}
\begin{picture}(4,1.7)
\put(0.5,1.4){$L\bc$}

\put(3.4, 1.4){$\bc$}

\put(0.5,0){$L\bc^{\ast}$}

\put(3.4, 0){$\bc^{\ast}$,}

\put(1.3, 1.55){\vector(1,0){1.75}}

\put(0.75, 1.2){\vector(0,-1){0.7}}

\put(3.5, 1.2){\vector(0,-1){0.7}}

\put(1.3, 0.15){\vector(1,0){1.75}}

\put(1.3, 0.3){\vector(2,1){1.8}}

\put(2, 0.25){$\varphi$}

\put(2, 1.65){$\phi$}

\put(2, 0.8){$\psi$}

\put(3.7, 0.8){$\exp$}

\put(0.1, 0.8){$\exp$}
\end{picture}
\end{equation}
where $\varphi \in \homoz$ is obtained by composition, and $\phi \in
\mathcal{A}$ is the Lie algebra homomorphism of $\psi$ and
$\varphi$. Let $L^0\bc^{\ast} \subset L\bc^{\ast}$ be the subgroup
of loops with winding number $0$. Given $\phi \in \mathcal{A}$, the
homomorphism $L^0\bc^{\ast} \ni x \mapsto \phi(\ln(x)) \in \bc$ is
independent of the choice of $\ln(x)$ and therefore well defined,
from which it follows that the space of $\phi$ (resp. $\varphi$)
generated as in (\ref{psi}) is all of $\mathcal{A}$ (resp.
$\homoz$). Thus $\homoz$, as the quotient space of $\homoc$ modulo a
discrete subgroup, is a complex Lie group.

Let $\psi \in \homoc$ and $g_{ab} \in G \subset LG$ such that
$g_{ab}(a)=\infty$ and $g_{ab}(b)=0$. Consider $\psi$ as a function
on $LU_{\infty 0}$ and define $\fh^{\psi}_{ab} =g_{ab}^{\ast}\psi
\in \co(LU_{ab})$. The property $\psi|_{\bc^{\ast}} \equiv 0$
implies that $\fh^{\psi}_{ab}$ is independent of the choice of
$g_{ab}$. In particular,
\begin{equation} \label{explicit}
\fh^{\psi}_{ab}(x)=\left\{ \begin{array}{ll}  \psi(x-b), &
a=\infty;\\ -\psi(x-a), & b=\infty;\\
\psi\left((x-b)(x-a)^{-1}\right), & \mathrm{otherwise}.
\end{array}\right.
\end{equation}
Thus $\fh^{\psi}=(\fh^{\psi}_{ab})$ is a \v{C}ech 1-cocycle of the
covering $\fu$. We claim that $\fh^{\psi}_{ab}$ is cuspidal with
respect to the constant $0$. It follows from the definition that
$\fh^{\psi}$ is $G$-fixed, so we only need to show the cuspidal
property of $\fh^{\psi}_{\infty 0}=\psi$. Indeed,
$$\lim_{\bc \ni \lambda \to \infty} \psi(x+\lambda)=\lim_{\lambda
\to \infty} \psi\left(x/\lambda+1\right)=0.$$ Therefore $\fh^{\psi}
\in \fH^{G}$.

If we replace $\psi$ by $\varphi \in \homoz$ in above, then the same
procedure yields a multiplicative \v{C}ech 1-cocycle $\fh^{\ast
\varphi}=(g_{ab}^{\ast} \varphi)$ $\in$ $(\fH^{\ast})^G$, which is
the same 1-cocycle as the one constructed in \cite[(2.4)]{z}.

\begin{lem}  \label{hom-cycle}
Let $\psi$, $\phi$ and $\varphi$ be as in (\ref{psi}).
\begin{enumerate}
\item[\textup{(a)}]  If we consider elements of $\mathcal{A}$
as functions on $\bc \times L\bc$ independent of the first variable,
then $\mathcal{A}=\ff^G$ and $\alpha(\fh^{\psi})=
\alpha^{\ast}(\fh^{\ast \varphi})=\phi$.

\item[\textup{(b)}] The map $\psi \mapsto \fh^{\psi}$ (resp.
$\varphi \mapsto \fh^{\ast \varphi}$) is an isomorphism $\homoc$
$\to$ $\fH^G$ (resp. $\homoz$ $\to$ $(\fH^{\ast})^G$) of complex
locally convex spaces (resp. complex Lie groups).
\end{enumerate}
\end{lem}
\begin{proof}
(a) It follows from the definition of $\ff$ as in (\ref{ff}) and of
the $G$-action on $\ff$ as in (\ref{gaction}) that $\mathcal{A}
\subset \ff^G$. Let $F \in \ff^G$. By considering the action of the
M\"obius transformations $\zeta \mapsto \zeta+\lambda$, $\lambda \in
\bc$, and $\zeta \mapsto \lambda \zeta$, $\lambda \in \bc^{\ast}$,
on $F$, we can conclude that $F$ is independent of the first
variable and linear with respect to the second variable. The third
property in the definition of $\ff$ implies that $F|_{\bc} \equiv
0$. Thus $F \in \mathcal{A}$ and $\ff^G \subset \mathcal{A}$.

Next we compute $\alpha(\fh^{\psi})$ by (\ref{f}) and
(\ref{explicit}). Choose $a=\infty$ in (\ref{f}), and we only need
to compute the derivative at $\zeta=0$ (as $\alpha(\fh^{\psi}) \in
\ff^G$ is independent of $\zeta$). So
\begin{eqnarray*}
\alpha(\fh^{\psi})(-1/x) &=& \left.
\frac{d}{d\zeta}\right|_{\zeta=0} \fh^{\psi}_{\infty \zeta}(x)
=\lim_{\zeta \to 0}\frac{\psi(x-\zeta)-\psi(x)}{\zeta}
\\ &=& \lim_{\zeta \to 0}
\frac{\psi\left(1-\zeta/x\right)}{\zeta}=\phi\left(-1/x\right),
\end{eqnarray*}
where $x \in LU_0$, and the second limit above represents the
derivative of $\psi$ at $1 \in L\bc^{\ast}$ in the direction of
$-1/x$. Thus $\alpha(\fh^{\psi})=\phi$. Note that
$\exp(\fh^{\psi})=\fh^{\ast \varphi}$. It follows from the
definition of $\alpha^{\ast}$ that $\alpha^{\ast}(\fh^{\ast
\varphi})=\phi$.

(b) As $\psi=\fh^{\psi}_{\infty 0}$, the map $\psi \mapsto
\fh^{\psi}$ is clearly injective. Let $\fh \in \fH^G$. By (a) there
exists $\fh^{\psi}$ such that $\alpha(\fh^{\psi})=\alpha(\fh) \in
\ff^G$. Thus $\fh -\fh^{\psi}=c (\mathrm{ind}_{ab})=\fh^{c \,
\mathrm{ind}_{\infty 0}}$, where $c \in \bc$ is a constant and
$\mathrm{ind}_{\infty 0} \in \homoc$. So $\fh=\fh^{\psi+c \,
\mathrm{ind}_{\infty 0}}$. Since $\fh^{\psi}_{ab}=g^{\ast}_{ab}
\psi$ and $\psi=\fh^{\psi}_{\infty 0}$, the isomorphism $\psi
\mapsto \fh^{\psi}$ is topological. The conclusion for the map
$\varphi \mapsto \fh^{\ast \varphi}$ follows from similar arguments.
\end{proof}

It follows from the relation between $\homoz$ and $(\fH^{\ast})^G$
and Proposition \ref{2t} that the quotient topology of $\homoz$
inherited from $\homoc$ is the same as the compact-open topology.

\vspace{3mm}

\noindent {\it Proof of Theorem \ref{main} \,} Since each bundle
$\Lambda \in \pic$ is a product of a $LG$-fixed bundle and a bundle
in $\picz$, see \cite[(2.5)]{z}, the conclusion $\pic^G=\pic^{LG}$
follows from $\picz^{G}=\picz^{LG}$. Let $\Lambda \in \picz^G$. By
Proposition \ref{iso} and Lemma \ref{hom-cycle}(b) $\Lambda$ is
generated by a 1-cocycle of the form $\fh^{\ast \varphi}$, $\varphi
\in \homoz$, the same one as in \cite[(2.4)]{z}. By
\cite[Propositions 2.2, 2.1]{z} $\Lambda$ is $LG$-fixed.

Let $[f] \in \dg^G$ and $[f_0]$ a generator of $\ker(\fe) \subset
\dg^{LG}$, see (\ref{seq}). Since $\fe([f])$ is $LG$-fixed, we have
$g[f]=[f]+m(g)[f_0]$, where $g \in LG$, and $m: LG \to \bz$ is a
group homomorphism, constant on each component of $LG$. Thus
$m\equiv 0$ and $[f] \in \dg^{LG}$.  \qed

\end{document}